\documentclass[12pt]{article}
\usepackage{amsmath}
\usepackage{amsfonts}
\usepackage{amssymb}
\usepackage{graphicx}

\numberwithin{equation}{section}
\newcommand{\R}{\mathbb R}

\begin{document}

\begin{center}
{\bf \Large{The Problem of Moments: A Bunch of Classical \\
 Results With Some Novelties}}
\end{center}

\vspace{0.7cm}
\noindent
{\sc Pier Luigi Novi Inverardi$^{a}$, \ Aldo Tagliani$^{a}$, Jordan M. Stoyanov$^{b}$}\footnote{E-mails: pierluigi.noviinverardi@unitn.it, aldo.tagliani@unitn.it, \ stoyanovj@gmail.com}

\vspace{0.2cm}
\noindent 
{\small $^{a}${\it Department of Economics $\&$ Management, University of Trento, 38100 Trento, Italy\\
$^{b}$Institute of Mathematics $\&$ Informatics, Bulgarian Academy of Sciences, 1113 Sofia, Bulgaria\\
and Faculty of Mathematical Sciences, Shandong University, Jinan 250100, P.R. China}}

\vspace{0.8cm}\noindent
{\bf \emph{Dedicated to Professor Christian Berg on the occasion of his $79\frac{1}{4}$ birthday! With great respect for his fundamental contribution to 
the  Moment Problem!}}

\vspace{0.4cm}\noindent
\begin{center}
{\bf Paper published in Symmetry (MDPI), 2023, 15, 1743, \ (19 pages)}\\
{\bf DOI: 10.3390/sym15091743}\\
{\bf Version of 01 October 2023} \ [corrected one typo and a few inaccuracies] 
\end{center}

\vspace{0.5cm}
\noindent 
{\small 
{\bf Abstract:}  We summarize  significant classical results on (in)determinacy of measures in terms of their finite positive integer order moments. 
 Well-known is the role of the smallest eigenvalues of Hankel matrices, 
starting from Hamburger's results a century ago and ending with the great progress made only in recent times by C. Berg and collaborators.  
  We describe here known results containing necessary and sufficient conditions for moment (in)determinacy in both Hamburger and Stieltjes moment problems. 
	In our exposition we follow an approach different from that commonly used.  There are  novelties well complementing the existing theory. 
	Among them are: (a) to emphasize on the geometric interpretation of the indeterminacy conditions;   (b) exploit fine properties of the eigenvalues 
	of perturbed symmetric matrices allowing to derive new lower bounds for the smallest eigenvalues of Hankel matrices; 
 these bounds are used for concluding indeterminacy; (c)   provide new arguments to confirm classical results; 
(d) give new numerical illustrations involving commonly used probability distributions.


  }

\vspace{0.1cm}\noindent
Received:  10 August 2023 \quad Revised: 31 August 2023  \quad Accepted: 01 September 2023

{\small 
\section{Introduction. Preliminaries} 
}

In order to make the paper self-contained and easy to follow, we first provide the basics. 
Within the classical moment problem, or  the problem of moments, we deal with two main questions whose answers are known and available in the literature.

\vspace{0.2cm}\noindent
\subsection{Basic terminology. Main questions} 

\vspace{0.2cm}\noindent
{\bf Question  1.}  (existence). Is there  
a bounded positive measure $\mu$ with a specified support ${\rm U} := {\rm supp}\,(\mu), \ {\rm U} \subset \R=(-\infty, \infty),$ such that a given 
infinite sequence of real numbers $\{m_k\}_{k=0}^{\infty}$ is the moment sequence of $\mu$, i.e.,  
 $m_k=m_k(\mu)=\int_{\rm U} x^k\,{\rm d}\mu(x)$ for 
all $k = 0, 1, 2, \ldots ?$ 

We assume, of course, that $\int_{\rm U} |x|^k\,{\rm d}\mu(x) < \infty,$ for all $ \ k=1, 2, \ldots;  m_k$ is the the $k$th order moment of $\mu.$ 
The answer to Question 1 is well-known; details given below.  

\vspace{0.2cm}\noindent
{\bf Question 2.} (uniqueness). Is $\mu$ the only measure with the moments $\{m_k\}_{k=0}^{\infty}$? 

\vspace{0.2cm}
If the answer to Question 2 is `yes', we say that $\mu$ 
is {\it moment determinate}, or, that $\mu$ is uniquely determined by its moments. Otherwise, if the answer is `no', $\mu$ is {\it moment indeterminate}, or, $\mu$ 
is not determined uniquely by its moments. 
In such a case, available is the following deep and non-trivial result, see Berg-Christensen (1981): 

\vspace{0.2cm}\noindent
{\bf General result.} {\it Suppose that $\mu$ is a measure with finite all moments of positive integer order. If $\mu$ is non-unique, 
then there are infinitely many measures of any kind, discrete, absolutely continuous or singular, 
all with the same moments as $\mu.$  }

\vspace{0.1cm}
There is a long, rich and amazingly interesting history which originates in works by P.L. Chebyshev (1832--1894) and A.A. Markov (1856--1922). 
The systematic development of the moment problem is due to 
T.J. Stieltjes (1856--1894), see his memoir Stieltjes (1894). In this remarkable work he was the first to show that the answer to Question 2  can be `no' 
by describing explicitly 
different measures on $\mathbb{R}_+$  sharing the same moments.    

 The answers to the above questions depend on both, the support ${\rm U}$ of  $\mu$ and the moment 
 sequence $\{m_k\}_{k=0}^{\infty}$. Adopted in the 
literature are the following names:   {\it Hausdorff moment problem}, if $\rm U$ is bounded;  
{\em Stieltjes moment problem}, if $\rm U$ is unbounded and ${\rm U} \subset \R_+$, and {\em Hamburger moment problem}, 
if $\rm U$ is unbounded and  ${\rm U} \subset \R$. 

If ${\rm U} = \R$, we say that $\{m_{k}\}_{k=0}^{\infty}$ is a {\it Hamburger moment sequence},  
 while for ${\rm U} = {\R}_+$,  $\{m_{k}\}_{k=0}^{\infty}$ is a {\it Stieltjes moment sequence}. 
We do not deal with the Hausdorff moment problem since any such a measure, when exists,  is uniquely determined by its moments.

Widely known references are the books by Shohat-Tamarkin (1943), Akhiezer (1965), Berg-Christinsen-Ressel (1984) and Schm\"udgen (2017); 
see also Simon (1998), Sodin (2019) and Olteanu (2023).  
These sources contain comprehensive details about a series of remarkable results paving 
the progress in  moment problems for more than a century.  

In this paper we use standard terminology and notations generally accepted in analysis in works on moment problems. 
Telling that  $\{m_{k}\}_{k=0}^{\infty}$ is a moment sequence always means that there is a measure $\mu$ `behind', 
i.e., there exists $\mu$ which produces these moments.

We need a few words for terminology clarity: If we are given a measure $\mu$ with finite moments, then $\mu$ produces its only one moment 
sequence $\{m_k\}_{k=0}^{\infty}$. 
 It is the measure $\mu$ which is either determinate or indeterminate, hence no reasons 
to stick `determinate' or `indeterminate' to the moment sequence.

 We are interested in the (in)determinacy property of a measure $\mu$ with unbounded support and finite all moments   
$\{m_{k}\}_{k=0}^{\infty}$. Notice, in general, a Stieltjes moment sequence $\{m_{k}\}_{k=0}^{\infty}$ can also be considered as a Hamburger moment sequence.
One important, not intuitive and nontrivial fact is the possibility for  a measure $\mu$ with finite moments $\{m_k\}$ to be  
determinate in Stieltjes sense and indeterminate in Hamburger sense. We refer, e.g., to Shohat-Tamarkin (1943), p. 75,  Berg-Valent (1994), p. 165, 
Schm\"udgen (2017), p. 183; see also Theorem 3.2 given below in Section 3.

\subsection{Hankel matrices and their smallest eigenvalues} 

For any moment sequence  $\{m_{k}\}_{k=0}^{\infty},$ we define a few  infinite sequences of {\it Hankel matrices}, 
namely,  $\{H_{n}\}_{n=1}^{\infty}$, called a `basic' Hankel matrix, and  $\{H_{n,p}\}_{n=1}^{\infty}$, called  a `$p$-shifted' 
Hankel matrix. Recall that $H_n$ and $H_{n,p}$ are $(n+1) \times (n+1)$ matrices defined as follows: 
\[
	H_{n}=(m_{i+j})_{i,j=0}^n \quad \mbox{ and } \quad  H_{n,p}=(m_{i+j+p})_{i,j=0}^n, \quad p=1, 2,  3, 4. 
\]

The basic Hankel matrix  $H_n$ is based on all moments $m_0, m_1, \ldots, $ while \ $H_{n,p},$ for $p=1, 2, 3, 4,$ 
the `shifted' Hankel matrices are formed as follows: $H_{n,1}=(m_{i+j+1})_{i,j=0}^n$ 
is based on the `shifted' moment sequence $\{m_1, m_2, \ldots \}$ which is generated by the measure $\mu_1$ with ${\rm d}\mu_1 = x{\rm d}\mu$;  
$H_{n,2}=(m_{i+j+2})_{i,j=0}^n$ is based on the `shifted' moment sequence $\{m_2, m_3, \ldots \}$ 
generated by the measure $\mu_2$ with ${\rm d}\mu_2 = x^2{\rm d}\mu$; similarly for $H_{n,3}$ and $H_{n,4}$. 
For the determinants of Hankel matrices we use  the following notations: 
\[
D_{n}:={\rm det}\,(H_{n}) \quad \mbox{ and } \quad D_{n,p}:= {\rm det}\,(H_{n,p}).
\]

\vspace{0.1cm}
In this paper we use $c_0, c_1, c_2$ to denote positive constants which depend on some 
fixed moments, but we omit them as 
 explicit arguments. For simplicity, 
if  the moments $m_{p+1}, \ldots, m_{p+2n}$ are  fixed and we allow $m_{p}$, the moment preceding $m_{p+1}$, to `vary', then 
instead of the full notations $H_{n,p}(m_p, m_{p+1}, \ldots, m_{p+2n})$ and 
$D_{n,p}(m_p, m_{p+1}, \ldots, m_{p+2n})$, we  write  $H_{n,p}(m_p)$  and $D_{n,p}(m_p).$

\smallskip
Recall a fact from Shohat-Tamarkin (1943), Theorems 1.2 - 1.3, which is related to our Question 1: 
A sequence of real numbers $\{m_{k}\}_{k=0}^{\infty}$ is the moment sequence of a measure  
 $\mu$ with support on $\R$ (Hamburger case) if and only if all Hankel matrices $\{H_n\}$ are non-negative definite, which is equivalent to the non-negativity of 
all determinants $\{D_n\}.$ If the support 
of $\mu$ is ${\R}_+$ (Stieltjes case), we have a similar statement, however now we need the two sequences of Hankel matrices, $\{H_n\}$ and $\{H_{n,1}\},$ 
to be non-negative definite, or, in terms of their determinants, that $D_n \geq 0$ and $D_{n,1}\ge0, n=0, 1, 2, \ldots$ 

In both cases, Hamburger and Stieltjes, if for some $n_0$ we have $D_{n_0}=0$, then $D_n=0$ for all $n>n_0$ and  $\mu$ is a discrete measure 
concentrated on a finite number of points. We also say that $\mu$ has a finite spectrum, i.e the support is bounded, which implies that the measure $\mu$ is unique. 

We are interested in the phenomenon `indeterminacy' (`nonuniqueness') of the measure $\mu$. This means that the support of $\mu$ has to be unbounded 
and its support not to be reducible to a finite set of points. This is why we deal with 
Hankel matrices which are strictly positive definite, equivalently, their determinants are strictly positive, $D_n>0, D_{n,1}>0, n = 0, 1, 2, \ldots.$

Since $D_{n}$, as well  $D_{n,1}$, are positive, then all eigenvalues of the Hankel matrices $H_n$ and $H_{n,p}$ are positive. 
We use the notation   $\lambda_k(H_{n,p}), \ k=1, 2, \ldots, n+1,$ for the $k$th largest eigenvalue of $H_{n,p}$. Thus
\begin{equation}\label{eigen}
 0 \leq \lambda_{min}(H_{n,p}) := \lambda_1(H_{n,p}) \leq \lambda_2(H_{n,p}) \leq \ldots \leq \lambda_{n+1}(H_{n,p}) := \lambda_{max}(H_{n,p}).
\notag
\end{equation}

It is well-known, see, e.g., Berg-Chen-Ismail (2002), that the smallest eigenvalue of the matrix $H_{n,p}$ is given by the Rayleigh relation 
\begin{equation}\label{Rayl}
 \lambda_{1}(H_{n,p})=\min \Big \{
 \Sigma_{j,k=0}^n\; m_{j+k+p}\;v_j v_k: \ v_0, v_1, \ldots, v_n \in \R, \ \Sigma_{i=0}^n v_i^2 = 1 \Big\}, 
 \notag
\end{equation}
and that the positive numerical sequence  $\{ \lambda_1(H_{n,p})\}_{n=1}^{\infty}$ is decreasing in $n$.

The smallest eigenvalues of Hankel matrices are fundamentally involved when studying (in)determinacy of measures; see the historical paper by 
Hamburger (1920) and the more recent work by Chen-Lawrence (1999), Berg-Chen-Ismail (2002) and Berg-Szwarc (2011).   

 To mention, in the case of indeterminacy the inverse Hankel matrices are involved.
Indeed, Berg-Chen-Ismail (2002) obtained the lower bound for the smallest eigenvalue in terms of the trace of the inverse matrix (see their eqns. 1.10, 1.11, 1.12), 
starting from the fact that $1/\lambda_{\min}$ of the Hankel matrix is equal to $\lambda_{\max}$ of the inverse of that Hankel matrix.

	\subsection{Two classical results} 
	
	For a very long time available in the literature were  classical results expressed in terms of the smallest eigenvalues of Hankel matrices; see,  
	e.g., Hamburger (1920).  		
	A remarkable progress was made only in more recent times by Berg-Thill (1991) and Berg-Chen-Ismail (2002). These authors proved fundamental results 
	which can be summarized as follows (the letter `H' stands for Hamburger, `S' stands for Stieltjes):  

\vspace{0.2cm}\noindent
{\bf Classical Result H.} 
{\it In the Hamburger moment 
problem, the measure $\mu$ is uniquely determined by its moments $\{m_k\}_{k=0}^{\infty}$ {\bf \emph{if and only if}} the sequence 
of the smallest eigenvalues of the basic Hankel matrices $\{H_n\}_{n=1}^{\infty}$ converges to zero as $n \to \infty:$ 
 $\lim_{n \to \infty} \lambda_1(H_n) = 0.$ 

Equivalently, $\mu$ is indeterminate by its moments $\{m_k\}_{k=0}^{\infty}$ {\bf \emph{if and only if}} the sequence 
of the smallest eigenvalues of the basic Hankel matrices $\{H_n\}_{n=1}^{\infty}$ converges to a strictly positive number as $n \to \infty:$ 
 $\lim_{n \to \infty} \lambda_1(H_n) = c_0,\ c_0 >0.$ 
 }

\vspace{0.2cm}\noindent
{\bf Classical Result S.} 
{\it In 
the Stieltjes moment problem, the measure $\mu$ is non-uniquely determined by its moments $\{m_k\}_{k=0}^{\infty}$ {\bf \emph{if and only if}} 
the sequences  
of the smallest eigenvalues of the basic Hankel matrices $\{H_n\}_{n=1}^{\infty}$ and of the shifted Hankel matrices $\{H_{n,1}\}_{n=1}^{\infty}$ 
both converge to strictly positive numbers:  $\lim_{n \to \infty} \lambda_1(H_n) = c_0,$ \ $\lim_{n \to \infty} \lambda_1(H_{n,1}) = c_1.$ 

Equivalently, the measure $\mu$ is determinate by its moments $\{m_k\}_{k=0}^{\infty}$ {\bf \emph{if and only if}} at least one of the 
sequences  
of the smallest eigenvalues of the basic Hankel matrices $\{H_n\}_{n=1}^{\infty}$ and of the shifted Hankel matrices $\{H_{n,1}\}_{n=1}^{\infty}$ 
converges to zero as $n \to \infty.$ }

Available in the literature are equivalent variations of the formulations of the above results. The proofs, however, may rely on different ideas and techniques.

\subsection{About the novelties in our approach}

Crucial in our approach is to exploit the following : \\
$\bullet$ \ The {\bf \emph{geometric interpretation}} of the indeterminacy conditions as developed by Merkes-Wetzel (1976).\\
$\bullet$ \ Properties of the {\bf \emph{eigenvalues of perturbed symmetric matrices}} in 
the spirit of Golub-Van Loan (1996) and Wilkinson (1995).

Both these are among the novelties in our exposition. They are properly used and combined with results from Shohat-Tamarkin (1943),  
Akhiezer (1965) and Schm\"udgen (2017) and  
a frequent referring  to Berg-Chen-Ismail (2002) or Berg-Thill (1991). Going this way we arrive at 
a unified presentation of classical results in both Hamburger and Stieltjes cases.  

As far as we are aware, there is no work, until now, giving such a presentation of most significant classical results on 
moment problems based on  ideas and techniques similar to those used in this paper.  
We found a little strange that the paper by Merkes-Wetzel (1976) 
was somehow neglected for a long time. It is not in the list of references in  
papers and books written by leading specialists on the moment problem.  
The only proper citation and comments are given by Wulfsohn (2006). 
In our opinion the 
geometric interpretation of the indeterminacy conditions has a value by its own, it is fresh and convincing, and deserves attention.  
The idea is quite simple. Based on the complete moment sequence $\{m_{k}\}_{k=0}^{\infty}$ we build up the so-called {\it parabolic limit region} in the plane,  
and then we look at the position of the point $(m_0, m_1)$. All depends on where this point is located: inside or outside of the region, or on its boundary.
  Later on we give details and clear graphical illustrations. 

We exploit intensively several properties of perturbed symmetric matrices, which allows to derive new lower bound used to conclude the indeterminacy property.  
Our bound is comparable with the lower bound derived in
 Berg-Chen-Ismail (2002) by using orthogonal polynomials. 
 
We provide a little different arguments, based on Krein-Nudelman (1977) for concluding the determinacy property.

\subsection{Moment determinacy in Probability theory} 

 It is worth mentioning that there are results  which are of the sort ``{\bf \emph{if and only if}}". 
Usually they are compactly formulated, fundamental in their content, and mathematically beautiful. However,  such results are difficult to prove 
and the conditions involved are practically impossible to check, hence the name  `uncheckable conditions'.  

If one assumes that $m_0 = 1$, i.e., that the total mass is $\mu(\rm U)=1$, then $\mu$ is a probability measure.
Well-known is the important role played by the moments in Probability and Statistics, and especially in their applications. 
This is why a special attention has been paid over a century on finding 
another sort of  `relatively easier' conditions which are only 
sufficient or only necessary for either determinacy or indeterminacy of a probability distribution. Nowadays, 
a variety of `checkable conditions' (Cram\'er, Hardy, Carleman, Krein) are available in the literature. 
The checkable conditions have their analytical value and are more than useful in several applied areas, see, e.g., Janssen-Mirbabayi-Zograf (2021). 

The paper by Lin (2017) is a rich and valuable source of information on 
classical and recent results on moment determinacy of probability distributions; see also Stoyanov-Lin-Kopanov (2020). 
The present paper is intrinsically related to another subsequent paper which is in preparation, see Lin-Stoyanov (2023).

\subsection{Structure of this paper} 

 The rest of the paper is organized as follows. In Section 2 we treat the Hamburger case and discuss conditions for   
(in)determinacy. Based on the geometric interpretation of the indeterminacy conditions, 
we re-derive in a different way already known results by Berg-Chen-Ismail (2002).  
  In Section 3 we follow the same line of reasoning  and establish results in the Stieltjes 
case announced in Berg-Thill (1991). In both cases we provide necessary and sufficient determinacy 
	conditions in terms of the asymptotic behavior, as $n \to \infty$, of two sequences of smallest eigenvalues, namely $\lambda_1(H_{n})$ 
	and  $\lambda_1(H_{n,1})$. We also provide a new lower bound for $\lambda_1(H_{n,1})$ which is related to the  
	indeterminacy of the measure involved. Section 4  presents details on  the smallest eigenvalues and their lower bounds 
	calculated in different ways. The numerical illustrations involve commonly used probability distributions. 
	Brief concluding comments  are given in Section 5.\\

\section{Hamburger moment problem} 
 
Before moving further, we discuss some known tools and results which will be used in the sequel. 
 
\subsection{ Limit parabolic region} 
 
From Shohat-Tamarkin (1943), p. 5 or Akhiezer (1965), p. 30, we know what $\{m_{k}\}_{k=0}^{\infty}$ is the moment sequence 
of a measure $\mu$ with  
support  $\R$ if and only if $D_n \ge 0$ for all $n=0,1,2, \ldots$. In such a case an H-sequence is called positive non-negative. 
We are going to deal with the interesting case of strictly positive H-sequences, in which case all $D_n > 0.$

Suppose now that $\{m_k\}_{k=0}^{\infty}$ is a moment sequence for which we `keep fixed' the moments $\{m_2, m_3, \ldots \},$
while we treat as `varying continuously' the moments $m_0$ and $m_1.$ If letting $m_0=x, \ m_1=y$,  the whole moment sequence can 
be written as $\{x, y, m_2, m_3, \ldots\}$. The numbers $x$ and $y$, i.e.,  
the moments preceding $m_2$, can not be arbitrary.  Moreover, if for $p=1, 2, 3, 4$, the moments $m_{p+1}, m_{p+2}, \ldots$ are 
fixed, there is always a range for the possible values of the `previous' moment $m_p$, the moment just before $m_{p+1}$. E.g., 
the following two-sided bound for $m_p$ holds:  
\[
a_{p,n}^- \leq m_p \leq (m_{p+1})^{p/(p+1)}.
\] 
Here the second relation is the Lyapunov's inequality. The lower bound $a_{p,n}^-$, i.e.,  
the smallest possible value of $m_p$, is the unique number such that 
\begin{equation}
	D_{n,p}(a_{p,n}^{-},m_{p+1}, \ldots, m_{p+2n}) =0.
\notag
\end{equation}

We now turn to an H-moment sequence and the geometric interpretation of the indeterminacy conditions for the corresponding 
measure. Following Merkes-Wetzel (1976), we fix $n$ and  consider the moments $(m_0, m_1, m_2, \ldots, m_{2n})$ . 
We keep `untouched' $m_2, \ldots, m_{2n}$, and assume that $m_0=x$ and $m_1=y$ are `varying continuously'. For each $n$ the relation 
 $D_{n}(x,y) \geq 0$ defines a closed convex region, ${\cal P}_n := \{(x,y) \in \R^2: D_{n}(x,y) \geq 0\}$, which is bounded by a 
proper parabola with horizontal 
axis and vertex in the right-half plane. Since for $n=1, 2, \ldots,$ the Hankel matrices $H_{n}(x,y)$ are positive semidefinite, the 
regions ${\cal P}_n$ are nested, ${\cal P}_n \subset {\cal P}_{n-1}$. Of interest is their intersection 
${\cal P} := \bigcap_{n=1}^\infty {\cal P}_n$, called a {\bf \emph{limit parabolic region}}. One possibility is that ${\cal P}$ is  a 
`proper' closed region 
in the right-half plane such that ${\cal P}$ is bounded by a proper parabola and  containing the initially given  
 moments $(m_0, m_1).$ The other possibility is ${\cal P}$ to degenerate to just a ray, as explained below.

\subsection{Determinacy criteria and their geometric meaning}  
 
In this subsection we provide a geometric interpretation of the determinacy criterion in the Hamburger case. The next two results, Theorems 1 and 2, play a fundamental role.

\smallskip
\noindent
{\bf Theorem 1.} (Shohat-Tamarkin (1943), Theorem 2.18).  
 	{\it Let $\mu$ be a measure associated with the strictly positive definite Hamburger moment sequence $\{m_k\}_{k=0}^{\infty}.$ 
	Then $\mu$ is determinate (unique) {\bf \emph{if and only if}} at least one of the two sequences of ratios  
\[
 	\frac{D_{n}}{D_{n-1, 2}} \quad \mbox{ and } \quad \frac{D_{n-1,2}}{D_{n-2, 4}}
\]
has a limit zero as $n \to \infty,$ i.e., either \ $\lim_{n \to \infty} \frac{D_n}{D_{n-1,2}}=0$ \ or \  
$\lim_{n \to \infty} \frac{D_{n-1,2}}{D_{n-2,4}}=0$.} 

\smallskip
\noindent
{\bf Theorem 2.} (Merkes-Wetzel (1976), Theorem 1). 
 	{\it Let  $\{x, y, m_{k}\}_{k=2}^{\infty}$ be a Hamburger moment sequence for the measure $\mu$. Then 
	 $\mu$ is indeterminate 
	{\bf \emph{if and only if}} the point $(x, y)$ is an interior point of the limit parabolic region $\cal P$.}\\

\smallskip
Theorem 1 has an interesting and enlightening geometric meaning. In fact, the following relations hold: 
\begin{equation}
\frac{D_{n}}{D_{n-1,2}} = m_0 - a_{0,n}^{-} \quad \mbox{ and } \quad \frac{D_{n-1,2}}{D_{n-2,4}}=m_2 - a_{2,n-1}^{-}. 
\notag
\end{equation} 
Here $a_{0,n}^-$ comes from the equation $D_{n}(a_{0,n}^{-}, m_1, \ldots, m_{2n})=0$, so the sequence $\{a_{0,n}^{-}\}_{n=1}^{\infty}$ 
is monotonic nondecreasing and, as $n \to \infty$, it is converging with $\lim_{n \to \infty} a_{0,n}^- :=a_{0,\infty}^{-} \leq m_0$. 
The number $a_{2,n-1}^{-}$ is defined from the relation $D_{n-1,2}(a_{2,n-1}^{-})=0$. It is remarkable that the difference  
$m_2 - a_{2,n-1}^{-}:=1/L_n$ is equal to the length of the `latus rectum' of the bounding parabola; see 
Merkes-Wetzel (1976), Lemma 1. The sequence of positive numbers $\{L_n\}_{n=1}^{\infty}$ is nondecreasing in $n$.   
Then, as $n \to \infty$ and $L_n \to \infty$, the length of the latus rectum $1/L_n$ tends to zero and the limit parabolic 
 region $\cal P$ becomes a ray; see Figure 1 (the red bold line). 

There are two possibilities for the measure $\mu$  with the moments $\{m_k\}$.

\noindent 
{\it Case 1.} (The measure $\mu$ is  H-\emph{indeterminate}). From Theorem 2 we have  that 
\[
 \lim_{n\to\infty}	\frac{D_{n-1,2}}{D_{n-2,4}} = m_2 - a_{0,\infty}^- := c_2 \quad
		{\hbox{and}} \quad \lim_{n \to \infty} \frac{D_{n}}{D_{n-1,2}} = 		m_0 - a_{0,\infty}^{-} :=c_0. 
\]
 It follows that the limit parabolic region is bounded by a non-degenerate parabola and, notice, the point $(m_0, m_1)$ is 
interior for the limit parabolic region $\cal P$, Figure 1. 
	
\begin{figure}[h] 
		\centering
		\includegraphics[width=0.8\linewidth]{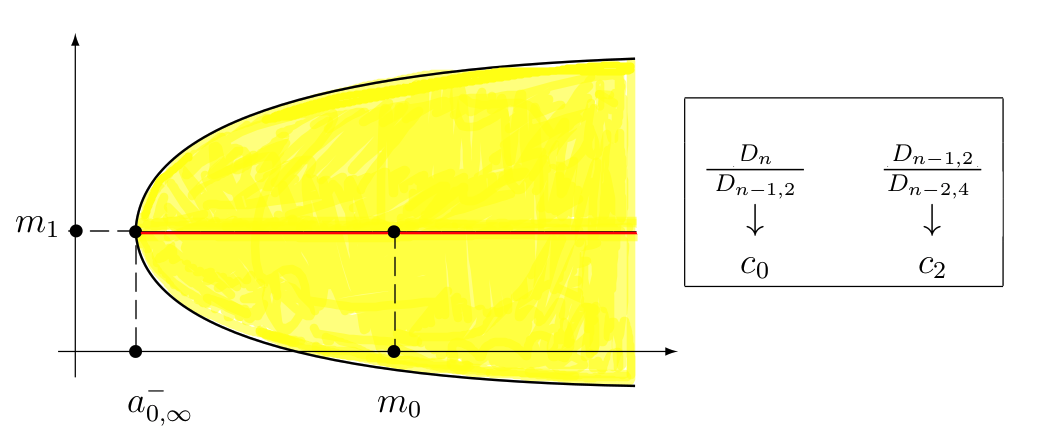}
		\caption{Parabolic region for H-indet case}
		\label{fig_1: H-indet}
	\end{figure}

\noindent 
{\it Case 2.} (The measure $\mu$ is  H-\emph{determinate}). 
    In view of Theorem 1, the limit parabolic region admits three distinct shapes. More precisely, Figure 2 
		is  analogues to Figure 1 		when  $c_0 = 0, \, c_2 > 0$. The two red bold lines in Figure 1 
		and in Figure 2 
		are referred to the two other (degenerate) parabolic regions (they become rays). This happens when either $c_0>0, \, c_2=0$ or 
		$c_0=0, \, c_2=0$, respectively.

\begin{figure}[h] 
    \centering
    \includegraphics[width=0.8\linewidth]{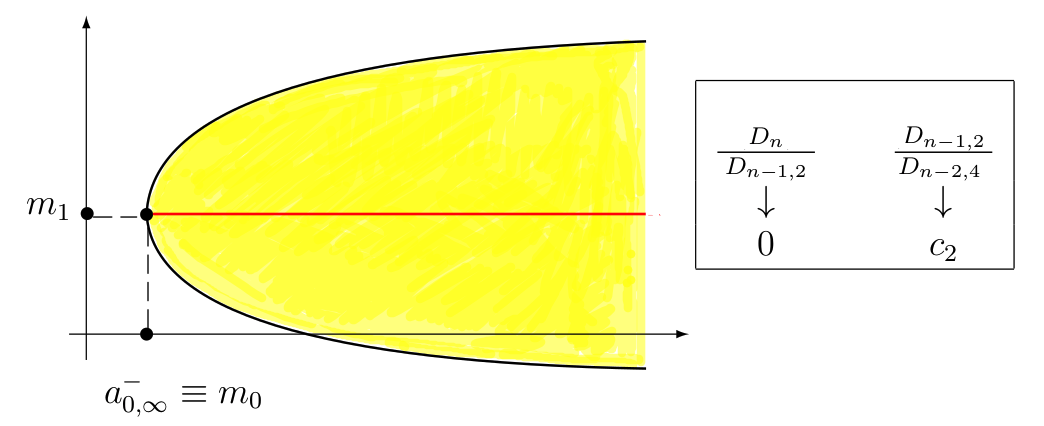}
		\caption{Parabolic region for H-det case}
		\label{fig_2: H-Det_(b)}
\end{figure}

   	\smallskip
   	Now, from Theorem 2  the following result promptly comes: Since the point $(m_0, m_1)$ is on the boundary of the limit parabolic 
		region, we have $\lim_{n \to \infty} D_n=0$ and hence $\lim_{n \to \infty} \lambda_1(H_n)=0.$ Consequently we arrive 
		at the following statement:
   	
   	\vspace{0.2cm}\noindent
		{\bf Corollary 1.} {\it  	
   		If the Hamburger moment sequence $\{m_k\}_{k=0}^{\infty}$ generates a measure $\mu$ which is determinate, then the monotonic 
			decreasing sequence $\{\lambda_1(H_n) \}_{n=1}^{\infty}$ has limit zero as $n \to \infty,$ i.e., $ \lim_{n \to \infty} \lambda_1(H_n)=0.$  
   	}

\subsection{Eigenvalues of perturbed symmetric matrices}  

Recall first that Berg-Chen-Ismail (2002)  use orthonormal polynomials to find a lower bound for the smallest 
eigenvalue $\lambda_1(H_{n})$, see their eq. (1.14). If taking large $n$ and  the fact that the limit $\lim_{n \to \infty} \lambda_1(H_{n})$
is strictly positive, they conclude the  indeterminacy of the  measure 
 associated with a Hamburger moment sequence.  

We are going now to derive properties of the smallest eigenvalues of a family of `perturbed' symmetric matrices. This allows 
to derive a lower bound for  $ \lim_{n \to \infty}\lambda_1(H_{n})$, which is a little different from the quantity 
found in the above cited paper. This, together with Theorems 1 and 2,   
leads to an alternative proof of Theorem 1.1 in Berg-Chen-Ismail (2002). Moreover, this approach can and will be followed  
to study the Stieltjes case; see Section 3.

 So, suppose that we have the H-moment sequence $\{m_k\}_{k=0}^{\infty}$ and  let $m_0>a_{0,n}^{-}$ for some fixed $n.$  
For a `small' number $\varepsilon \geq 0$ 
and ${\bf e}_1 =(1, 0, \ldots, 0)$ being the unit fundamental row-vector in $\R^{n+1}$ we define the following $(n+1) \times (n+1)$ matrix: 
\[
E \ = \ \varepsilon \,{\bf e}_1'\,{\bf e}_1, 
\]
which is a symmetric, non-negative matrix of rank 1. This $E$ will play the role of a `perturbation matrix'  
of a Hankel matrix. Consider now the $(n+1) \times (n+1)$ matrices 
\[
H_{n}(m_0) \quad \mbox{ and } \quad  H_{n}(m_0+\varepsilon)= H_{n}(m_0)+E.
\]
Here $H_{n}(m_0)$ is the Hankel matrix based on the moments $m_0, m_1, \ldots, m_{2n}$, while the matrix $H_{n}(m_0+\varepsilon)$ 
is based on $m_0 + \varepsilon, m_1, \ldots, m_{2n}$.

As a consequence of the simple form of $E$, the eigenvalue $\lambda_k(H_{n}(m_0))$ of the matrix $H_{n}(m_0)$ and 
the eigenvalue $\lambda_k(H_{n}(m_0  +\varepsilon))$ of the perturbed matrix 
  $H_{n}(m_0 +\varepsilon)$ are closely related to one other as follows (see, e.g., Wilkinson (1965), formula (41.8), p. 94-98 and 
	Golub-Van Loan (1996), Theorem 8.1.8, p. 397): 
\begin{equation}\label{Wilkinson}
    \lambda_k(H_{n}(m_0+\varepsilon))= \lambda_k(H_{n}(m_0)) +b_k\varepsilon  \quad \ \mbox{ for } \quad k=1, \ldots, n+1,  
\end{equation}
where 
\begin{equation}\label{Wilkinson1}
    \lambda_k(H_{n}(m_0)) \leq  \lambda_k(H_{n}(m_0+\varepsilon)) \leq  \lambda_{k+1}(H_{n}(m_0)), \quad k=1,...,n. 
\end{equation}   
It remains to tell that $b_1, \ldots, b_{n+1}$ are real numbers, each in the interval $[0,1],$ with sum   $\sum_{k=1}^{n+1} b_k=1$. 
The latter equality comes from \eqref{Wilkinson} and the obvious fact for traces of matrices: 
$ {\rm Tr}\,(H_{n}(m_0+\varepsilon))={\rm Tr}\,(H_{n}(m_0))+\varepsilon$.

Now we suitably specify \eqref{Wilkinson} by using the assumption that $m_0 > a_{0,n}^{-}$ and by setting $\varepsilon:=m_0 - a_{0,n}^{-}$. Hence 
 \eqref{Wilkinson} and \eqref{Wilkinson1} become
\begin{equation}\label{Wilkinson2}  
\lambda_k(H_{n}(m_0))= \lambda_k(H_{n}(a_{0,n}^-))+ b_k\,(m_0 - a_{0,n}^{-})
 \end{equation}
and
 \begin{equation}\label{Wilk1}  
 \lambda_{k}(H_{n}(a_{0,n}^{-})) \leq  \lambda_k(H_{n}(m_0)) \leq  \lambda_{k+1}(H_{n}(a_{0,n}^{-})), \quad k=1,...,n. 
 \end{equation}
 
 Applying Lagrange's theorem to the right-hand side of \eqref{Wilkinson2} on the interval $[a_{0,n}^{-}, m_0]$ we obtain 
 that for some number $\theta_n(m_0) \in (a_{0,n}^-, m_0)$ one holds 
\begin{equation}\label{Wilk}
\lambda_k(H_{n}(m_0))= \lambda_k(H_{n}(a_{0,n}^-))  +b_k(\theta_n(m_0);n)(m_0-a_{0,n}^{-}). 
\end{equation}

Eq. \eqref{Wilk} compared with \eqref{Wilkinson2} yields 
\[
b_k=\frac{\rm d}{{\rm d}x}\lambda_k(H_{n}(x))|_{x={\theta_n(m_0)}} = :b_k(\theta_n(m_0);n). 
\]  

 Let us list some consequences from \eqref{Wilk}. 

First, all  eigenvalues $\lambda_k(H_{n}(m_0))$
 are monotonically increasing with respect to $m_0$;
 since $0\leq b_k(\theta_n(m_0);n)\leq 1$, \eqref{Wilkinson} implies that  
	$\lambda_k(H_{n}(m_0))<m_0-a_{0,n}^{-}$, which means that each $\lambda_k(H_{n}(m_0))$, $k=1,...,n$, is shifted 
	by an amount which lies between zero and the positive number  $\varepsilon = m_0-a_{0,n}^{-}$. 
	
	Second, combining  \eqref{Wilk}  and \eqref{Wilk1} and taking $m_0 \to \infty$, for $k=1, \ldots, n$, 
	we see that $\lambda_k(H_{n}(m_0))$ takes a positive constant value, $\lambda_k(H_{n}(m_0)) = {\tilde c}_k>0$, while 
     $b_{n+1}(m_0; n) \to 1$. Hence, if $n \to \infty$ and $m_0 \to \infty$, we find that $b_k(m_0;n) \to 0$ for  $k=1, \ldots, n$ 
		which exhibits a different limit behavior compared with  $b_{n+1}(m_0; n) \to 1$. 
    Thus shows that  $\lambda_{n+1}(H_{n}(m_0))$ is asymptotically linearly increasing in $m_0$.      
     
   Third,  from the relations $D_{n}(m_0) = (m_0-a_{0,n}^{-})\,D_{n-1,2}$, coming from Theorem 1 and their
	geometric meaning, one has 
   \begin{equation}
   \ln D_{n}(m_0)=\ln \prod_{k=1}^{n+1} \lambda_k(H_{n}(m_0))=
   \sum_{k=1}^{n+1} \ln \lambda_k(H_{n}(m_0)).
   \notag
   \end{equation}
   Differentiating both sides with respect to $m_0$, we find 
   \begin{equation}\label{Weyl1}
   \frac{1}{m_0-a_{0,n}^{-}} =
   \sum_{k=1}^{n+1} \frac{1}{\lambda_k(H_{n}(m_0))} 
   \frac{\rm d}{{\rm d}m_0}\lambda_k(H_{n}(m_0))= 
   \sum_{k=1}^{n+1} \frac{b_k(m_0;n)}{\lambda_k(H_{n}(m_0))} 
   > \frac{b_1(m_0;n)}{\lambda_1(H_{n}(m_0))}
   \notag
   \end{equation}
   \noindent
   from which  
   \begin{equation}\label{Weyl2}
    b_1(m_0;n)(m_0-a_{0,n}^{-} )<\lambda_1(H_{n}(m_0)). 
   \end{equation}     
Combining \eqref{Weyl2} with \eqref{Wilk} shows that one holds $b_1(\theta_n(m_0);n)> b_1(m_0;n)$. Hence we conclude that the function
$ b_1(m_0;n)$ is monotonic decreasing as $m_0$ increases, which means that        
$\lambda_1(H_{n}(m_0))$ is a concave function.

Summarizing our findings above show that in this H-indeterminate case 
we deal with two  strictly positive quantities, namely 
 $m_0 - a_{0,\infty}^{-}>0$ and  $\lambda_1(H_{n})>0$. We use these facts to prove the following statement. 
 
\vspace{0.2cm}\noindent
{\bf Lemma 1.} {\it Assume that $m_0 > a_{0,\infty}^{-}$. Then, as $n \to \infty$, the sequence $\{b_1(m_0;n) \}_{n=1}^{\infty}$ admits a positive limit, 
denoted by $b_1(m_0;\infty)$, in the sense that for a suitable number $\tilde m_0>m_0,$ we will have 
\[
b_1(m_0;\infty) = \lim_{n \to \infty}\frac{\lambda_1(H_{n}(\tilde m_0))}{\tilde m_0-a_{0,n}^{-} }.
\]
}
{\bf Proof.}  We observe first that  $b_1(m_0;\infty)$ cannot be directly drawn by taking the limit in \eqref{Weyl2}, 
although each of the quantities $(m_0 - a_{0,n}^{-} )$ and $\lambda_1(H_{n})$ has a limit. Our arguments are given in the next 
three steps. 

\noindent
{\it Step 1.}  Consider \eqref{Wilk}. 
As $n\to\infty$, $\lambda_1(H_{n}(m_0))$ has a limit. Also, we have that  $m_0-a_{0,n}^{-}\to m_0-a_{0,\infty}^{-}$.
	As a consequence we have the relations
\[
	\lim_{n \to \infty} b_1(\theta_n(m_0);n) =: b_1(\theta_\infty(m_0);\infty)= 
	\lim_{n \to \infty}\frac{\lambda_1(H_{n}(m_0)) }{m_0-a_{0,n}^{-} }  
\]

Equivalently, both limiting quantities $\theta_\infty(m_0)$ and $b_1(\theta_\infty(m_0);\infty)$ exist, 
and moreover, we have the relations $a_{0,\infty}^{-} < \theta_\infty(m_0)<m_0$. \\
{\it Step 2.}  For fixed $n$, consider the interval $[a_{0,n}^{-}, m_0]$ and the function $\theta_n(m_0)$.
	From the concavity of $\lambda_1(H_{n}(m_0))$, as $m_0$ increases, each of the quantities $b_1(\theta_n(m_0);n)$ 
	$b_1(\theta_\infty(m_0);\infty)$ is decreasing. This implies that  
	both $\theta_n(m_0)$ and  $\theta_\infty(m_0)$ are increasing in $m_0.$ \\
{\it Step 3.} 
	Combining Step 1 and Step 2 we see that there exists a number $\tilde m_0>m_0$ such that $\theta_\infty(\tilde m_0)= m_0$.
	
 From \eqref{Wilk} with $m_0=\tilde m_0$ then it follows that
 \begin{equation}
  \lim_{n \to \infty} b_1(\theta_n(\tilde m_0);n)= b_1(\theta_\infty(\tilde m_0);\infty)=b_1(m_0;\infty)=\lim_{n \to \infty}
  \frac{ \lambda_1(H_{n}(\tilde m_0)) }{ \tilde m_0-a_{0,n}^{-} }.
  \notag
  \end{equation}
  This means that the number $b_1(m_0;\infty)$ indeed exists. $\Box$

 \smallskip
   From \eqref{Weyl2} we find the required lower bound for the eigenvalue $\lambda_1(H_{n})$:  
\begin{equation}\label{lamda1A}
   		b_1(m_0;\infty)(m_0 - a_{0,\infty}^{-})\leq \lim_{n \to \infty} \lambda_1(H_{n}),
\end{equation}  
which was one of our goals. 
With \eqref{lamda1A} 
in hands we are in a position to formulate and confirm the validity of the following well-known result, 
Theorem 1.1 in Berg-Chen-Ismail (2002). Notice, we are arriving at this result in a different way. 

\smallskip\noindent 
{\bf Theorem 3.} {\it Suppose that $\mu$ is a measure associated with the Hamburger moment sequence $\{m_{k}\}_{k=0}^{\infty}$  
and let $\lambda_1(H_n)$ be the smallest eigenvalue of the Hankel matrix $H_n.$ Then $\mu$ is determinate {\bf \emph{if and only if}}
    \begin{equation}\label{Hbound2}
    \lim_{n \to \infty} \lambda_1(H_{n}) = 0. 
     \end{equation}
		}
\noindent
{\bf Proof}. In one direction, if we assume that $\mu$ is H-determinate, 
	 from Corollary 1 we have $\lim_{n \to \infty}\lambda_1(H_n)=0$. In the other direction, if 
	$\mu$ is H-indeterminate, in view of Lemma 1, $\lambda_1(H_n)$ would have a positive lower bound \eqref{lamda1A}.  $\Box$
  
\bigskip
Most important, this is 
a {\bf \em{criterion for determinacy in the Hamburger moment problem}}.
Remarkable is the fact that involved are only the basic Hankel matrices $\{H_n\}_{n=1}^{\infty}$. 
Recall that the classical result of Hamburger (1920) involves the basic matrices $\{H_n\}_{n=1}^{\infty}$ {\bf and} the 
  shifted matrices $\{H_{n,2}\}_{n=1}^{\infty}$. The proof in Berg-Chen-Ismail (2002), as well our proof above, 
	do not rely on results in Hamburger (1920). 

We can say a little more. Recall the following known chain of relations from Akhiezer (1965), pp. 84--85: 
\begin{equation}
\lim_{n \to \infty}\frac{D_{n}}{D_{n-1,2}} = m_0 - a_{0,\infty}^{-}=\rho(0)= \frac{1}{\sum_{n=0}^\infty p_n^2(0)}, 
\notag
\end{equation}
where $p_n(x)$ is the $n$th orthonormal polynomial for the measure $\mu$. Then, if the function 
$\rho(x)= 1/\sum_{n=0}^\infty p_n^2(x) \not= 0$ for all $x \in \R$,  the measure $\mu$ is indeterminate. 
It remains to see that the lower bound in \eqref{lamda1A} is then replaced by  
\begin{equation}
\frac{ b_1(m_0;\infty)}{\sum_{n=0}^\infty p_n^2(0)}.
\notag
\end{equation}
In a sense, this lower bound is similar to the bound found by Berg-Chen-Ismail (2002), eqs. (1.14)--(1.15). In both cases 
the conclusion drawn is of course the same. 
  
\section{Stieltjes moment problem}   

In this case we develop a procedure which is similar to that followed in the Hamburger case, with some specifics. 
 We rely essentially on two known results, Theorems 4 and 5.

\smallskip
\noindent
{\bf Theorem 4.} (Merkes-Wetzel (1976), Lemma 3).  
{\it	The measure $\mu$ associated with the strictly positive definite Stieltjes moments sequence $\{m_k\}_{k=0}^{\infty}$ is determinate 
{\bf \emph{if and only if}} at least one of the sequences 
\[
\frac{D_{n}}{D_{n-1,2}}  \quad \mbox { and } \quad \frac{D_{n,1}}{D_{n-1,3}}
\]
	has a limit zero as $n \to \infty$, i.e.,  $\lim_{n \to \infty} \frac{D_n}{D_{n-1,2}}=0, \ \mbox{ or } \ 
	\lim_{n \to \infty} \frac{D_{n,1}}{D_{n-1,3}}=0.$}

\smallskip
It is useful to mention that Theorem 4 has the following geometric meaning: 
$\frac{D_{n,1}}{D_{n-1,3}} = m_1 - a_{1,n}^{-}$, where the number $a_{1,n}^{-}$ is the unique solution of the equation
 $D_{n,1}(a_{1,n}^{-})=0$. The numerical sequence $\{a_{1,n}^{-}\}_{n=1}^{\infty}$ is monotonic nondecreasing and, as $n \to \infty$,  
convergent to a limit, say $a_{1,\infty}^{-}$, where $a_{1,\infty}^{-} \leq m_1$.

\smallskip
\noindent
{\bf Theorem 5.} (Merkes-Wetzel (1976), Theorem 2) \label{Thm_5} 
	{\it The strictly positive definite Stieltjes moment sequence $\{m_k\}_{k=0}^{\infty}$ generates an indeterminate measure $\mu$  
	{\bf \emph{if and only if}} two conditions are satisfied: (i) the point $(m_0, m_1)$ is interior for the limit parabolic region $\cal P$; 	
	(ii) $a_{1,\infty}^{-} < m_1$.}

\smallskip
We recall that every Stieltjes sequence can also be considered as a Hamburger sequence; see, e.g., Chihara (1968).  Hence, the existence of 
the limit parabolic region $\cal P$ 
 is assured and defined by the two relations, $D_n(x,y) \geq 0$ and $D_{n,1}(y) \geq 0$. Her we use the notations $x=m_0$ and  $y=m_1$.  
Either ${\cal P}=\{(x,y): x \geq a_{0,\infty}^{-}, \ y = a_{1,\infty}^{-}\}$ is a ray, or $\cal P$ is the intersection of  proper limit 
parabolic regions  in the half plain $\{y \geq a_{1,\infty}^{-}\}$.

\smallskip
Consider the shifted Hankel matrix $H_{n,1}(m_1)$, its smallest eigenvalue $\lambda_1(H_{n,1}(m_1))$ and the perturbation matrix $E$, the same 
as previously defined in the Hamburger case. First, we want to show that in the S-indeterminate case the estimate of $\lambda_1(H_{n}(m_0))$ 
as $m_0$ varies and the estimate of $\lambda_1(H_{n,1}(m_1))$ as $m_1$ varies are equivalent procedures,    just  
 replace $D_{n}$ and $m_0$ with  $D_{n,1}$ and $m_1$. 

Now we need a result, which is a criterion for S-indeterminacy:  

\smallskip
{\it Suppose that $\mu$ and $\mu_1$ are measures associated with the moment sequence $\{m_{k}\}_{k=0}^{\infty}$ and the 
shifted moment sequence $\{m_{k+1}\}_{k=0}^{\infty}$, respectively. 
Then $\mu$ is S-indeterminate {\bf \emph{if and only if}}  both $\mu$ and $\mu_1$ are H-indeterminate.}

\smallskip
This statement is from Krein-Nudelman (1977), p. 199, P.6.8., where it is left as an exercise to the readers. For the sake of completeness 
we include here the proof. 

Indeed, assume that $\mu$ is S-indeterminate. This implies the H-indeterminacy of $\mu$ and S-indeterminacy of $\mu_1$. 
The latter yields H-indeterminacy of $\mu_1.$ If assuming that both $\mu$ and $\mu_1$ are H-indeterminate, we use
 Theorem 1 (formulated for determinate measures). Thus we have two limiting relations,  
		$\frac{D_{n}}{D_{n-1,2}} \to c_0$ and $\frac{D_{n,1}}{D_{n-1,3}} \to c_1$, where $c_0>0$ and $c_1>0.$  
		By Theorem 4 we conclude that $\mu$ is S-indeterminate.  
 
\smallskip
We can use Theorem 2 and describe alternatively the S-indeterminacy and also the H-indeterminacy in geometric terms. 
For this purpose we introduce two limit parabolic regions, ${\cal P}_{\rm H}$ in the Hamburger case, and ${\cal P}_{\rm S}$ in the 
Stieltjes case. With the convention $x=m_0, y=m_1, z=m_2$, we define: 
\[
{\cal P}_{\rm H} = \cap_{n=1}^\infty {\cal P}_n, \mbox{ where } {\cal P}_n = \{(x,y) \in {\mathbb R}^2: D_{n}(x,y) \geq 0\}, 
\]
\[
{\cal P}_{\rm S} = \cap_{n=1}^\infty {\tilde {\cal P}}_n, \mbox{ where } {\tilde {\cal P}}_n = \{(y,z) \in {\mathbb R}^2: D_{n,1}(y,z) \geq 0\}.
\]

As before, $\mu$ is a measure corresponding to the Stieltjes moment sequences $\{m_{k}\}_{k=0}^{\infty}.$ 
We have the following transparent interpretation:  

{\it The measure $\mu$ is {\rm S}-indeterminate, hence also {\rm H}-indeterminate, {\bf \emph{if and only if}} two conditions are satisfied: 
(i) the point $(m_0,m_1)$ is interior for the region ${\cal P}_{\rm H}$; (ii) the point $(m_1,m_2)$ is interior for the 
region ${\cal P}_{\rm S}.$ }

\smallskip
Since $a_{1,\infty}^{-}<m_1$, we refer to \eqref{lamda1A} and write down the following lower bound of the smallest 
eigenvalue $\lambda_1(H_{n,1})$ of the shifted Hankel matrix $H_{n,1}$:   
 \begin{equation}\label{Sbound2}
  	 b_1(m_1;\infty)\cdot(m_1 - a_{1,\infty}^{-}) \leq \lim_{n \to \infty}\lambda_1(H_{n,1}).
\end{equation}
Notice, this bound is related to the determinacy of the measure $\mu_1$ with  $ {\rm d}\mu_1 = x\,{\rm d}\mu.$

\smallskip 
 Let us  summarize the above findings: If a measure is S-indeterminate, it is also H-indeterminate. However, an S-determinate measure 
can be either H-determinate or H-indeterminate. Thus we have the cases, briefly discussed bellow.

\smallskip
\noindent
{\it Case 1.} ($\mu$ is S-indeterminate and H-indeterminate). From Theorem 1 and Theorem 4 
we have the inequalities $a_{0,\infty}^{-}< m_0$ and $a_{1,\infty}^{-} < m_1$. Then from \eqref{Hbound2} and \eqref{Sbound2} 
it follows that 
$\lim_{n \to \infty}\lambda_1(H_{n})=c_0>0$ and $\lim_{n \to \infty}\lambda_1(H_{n,1})=c_1>0$. 
Conversely, if $\lim_{n \to \infty}\lambda_1(H_{n})=c_0>0$ and $\lim_{n \to \infty}\lambda_1(H_{n,1})=c_1>0$, then the following 
two relations hold:  $a_{0,\infty}^{-}<m_0$ and $a_{1,\infty}^{-}<m_1$.  

\smallskip
\noindent
{\it Case 2.} ($\mu$ is S-determinate and H-indeterminate; see Merkes-Wetzel (1976), Corollary, p. 417). 
Since $\mu$ is S-determinate, Theorem 4 implies that $a_{1,\infty}^{-}=m_1$ and 
then that $\lim_{n \to \infty}\lambda_1(H_{n,1})=0$. 
Conversely, starting with $\lim_{n \to \infty}\lambda_1(H_{n,1})=0$, the relation $a_{1,\infty}^{-}=m_1$ 
follows from the S-determinacy. See Figure 3. 
      
 	\begin{figure}[h] 
		\centering
		\includegraphics[width=0.8\linewidth]{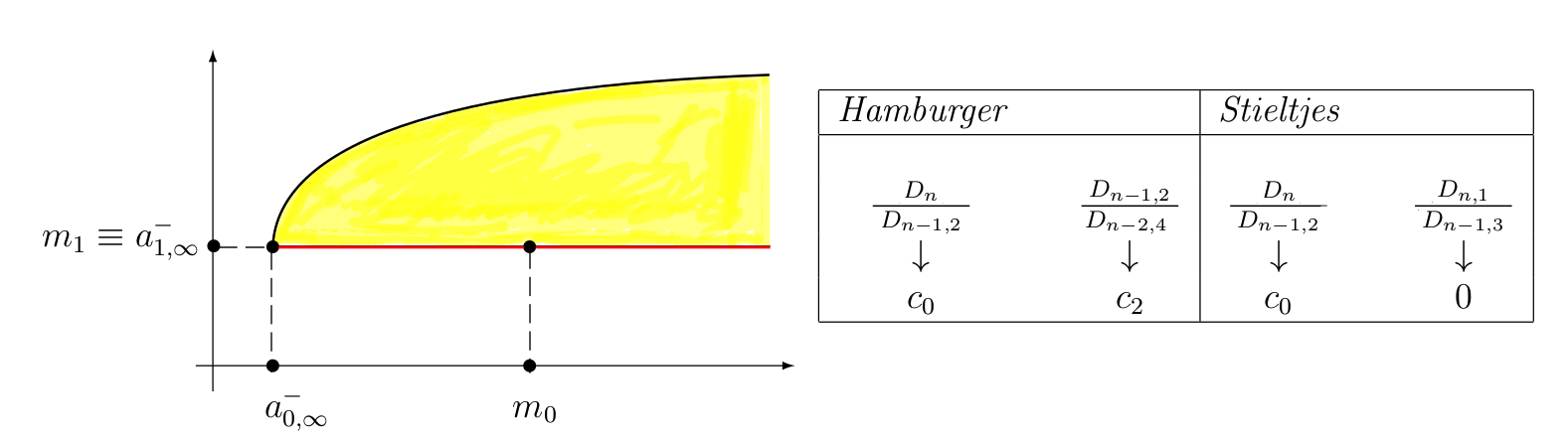}
		\caption{Parabolic region for S-det - H-indet case}
		\label{fig_4: H-Indet_S-Det_New}
	\end{figure}

\noindent 
{\bf Remark 1.} The S-determinate measure $\mu$ on $[0,\infty)$ in Case 2  is the Nevanlinna-extremal measure and the corresponding
Pick function coincides with a constant equal to zero. Hence $\mu$  is a discrete measure concentrated on the zeros of 
the $D$-function in the Nevanlinna parametrization, $D(x)=x\sum_{n=0}^\infty p_n(0)p_n(x)$. 
In particular, $\mu$ has a mass at 0. For details see Berg-Valent (1994), Remark 2.2.2, p. 178.

One possibility to construct a measure which is S-determinate and H-indeterminate is to start with a moment 
sequence $\{m_{k}\}_{k=0}^{\infty}$ associated with S-indeterminate measure, hence this moment sequence 
corresponds also to H-indeterminate measure.
Then the idea is to modify this sequence and get another one, $\{{\tilde m}_{k}\}_{k=0}^{\infty}$, 
associated with an S-determinate measure which is H-indeterminate. Such a specific construction is given in 
Sch\"umdgen (2017), Example 8.11, p. 183. It is shown (we do not give details here) how to calculate a 
proper constant $u>0$ and define the new moments by 
$\tilde m_k= \sum_{j=0}^k{k\choose j}(-u)^j m_{k-j}$ for $ k=0, 1, 2, \ldots. $ 
 Related relevant details can be found in Simon (1998), p. 96, Theorem 3.3.

\smallskip
\noindent
{\it Case 3.} ($\mu$ is S-determinate and H-determinate). For a given moment sequence, S-determinacy means that there is 
only one measure with support $[0,\infty)$. Regarding H-determinacy,  Corollary 1 provides an exhaustive answer. 

\smallskip
Clearly, important is the value of the limit   $\lim_{n \to \infty}\lambda_1(H_{n,1})$.  
Combining Theorems 1, 2, 3 and 4, we find that there are four possible limit parabolic regions, (a) - (d), 
and they are all feasible. Below are the details. 

\smallskip
\noindent 
(a) \ From Case 1 (S-indeterminate and H-indeterminate), we deal with a moment sequence $\{m_{k}\}_{k=0}^{\infty}$ 
whose associated discrete measure, say $\nu$, has a mass at 0 (this comes from the H-indeterminacy condition  
$a_{0,\infty}^{-} < m_0$). Next, consider the measure  $\tilde\nu$ 
related to $\nu$ via the relation  $\tilde\nu=\nu -\nu(\{0\})\delta_0$. The moments sequence $\{\tilde m_{k}\}_{k=0}^{\infty}$ 
of the measure $\tilde\nu$ differs from $\{m_{k}\}_{k=0}^{\infty}$ only at the zero-th entry (see Berg-Christensen, (1981), 
Theorem 7, p. 111). Hence $\tilde\nu$ has mass zero at 0, so that $\tilde a_{0,\infty}^{-}=\tilde m_0$. From Theorems 1 and 3  
and their geometric meaning we conclude for $\tilde\nu$ both properties, H-determinacy and S-determinacy. The limit parabolic region is 
nondegenerate with $(\tilde m_0,\tilde m_1)$ on its boundary and $\tilde a_{1,\infty}^{-}<\tilde m_1$, so that one 
holds $\lim_{n\to\infty}\lambda_1(H_{n,1})=c_1>0.$  See Figure 4. 

\begin{figure}[h] 
		\centering
		\includegraphics[width=0.8\linewidth]{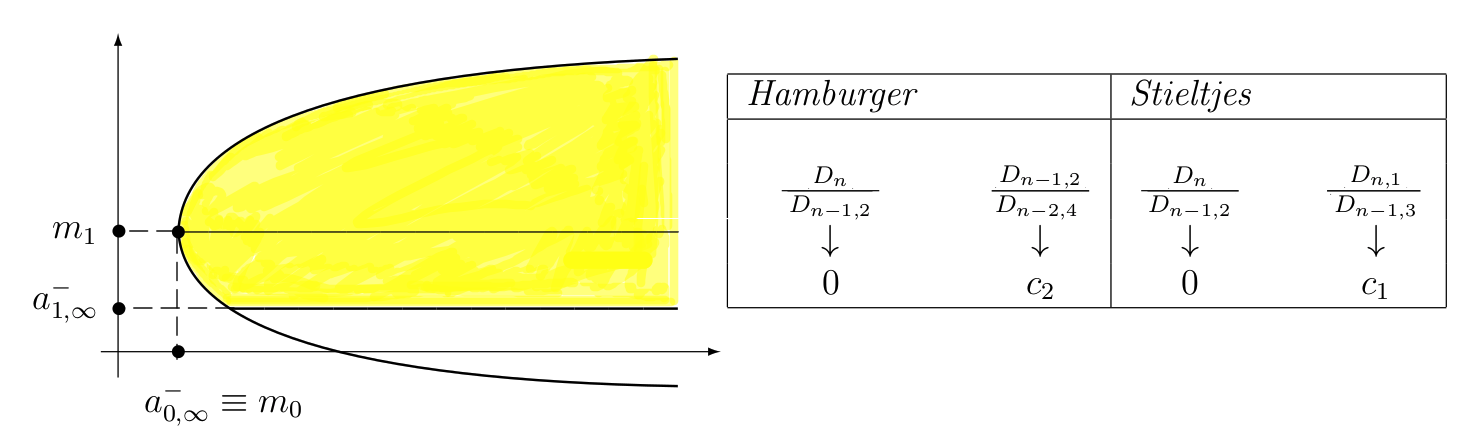}
		\caption{Parabolic region for S-det with H-det - (a)}
		\label{fig_5: H-Det-S-Det_New (a)}
	\end{figure}

\smallskip
\noindent 
(b) \ Here we have a nondegenerate limit parabolic region such that the point $(m_0, m_1)$ is on its boundary 
(which implies H-determinacy) and that $a_{1,\infty}^{-}= m_1.$ The latter implies S-determinacy and then 
$\lim_{n \to \infty} \lambda_1(H_{n,1})=0$.  
	 In this case the unique measure, say  $\tilde\mu,$ is related to the measure $\mu$, involved in the above Case 2 
	(S-determinate 	and H-indeterminate) by  the relation $\tilde\mu=\mu -\mu(\{0\})\delta_0$. By analogy with the 
	previous item (a), the measure $\tilde\mu$ has a moment  sequence $\{\tilde m_{k}\}_{k=0}^{\infty}$ which 
	 differs from $\{m_{k}\}_{k=0}^{\infty}$ only by the very first entry indexed by 0 (zero).
 
\begin{figure}[h] 
		\centering
		\includegraphics[width=0.8\linewidth]{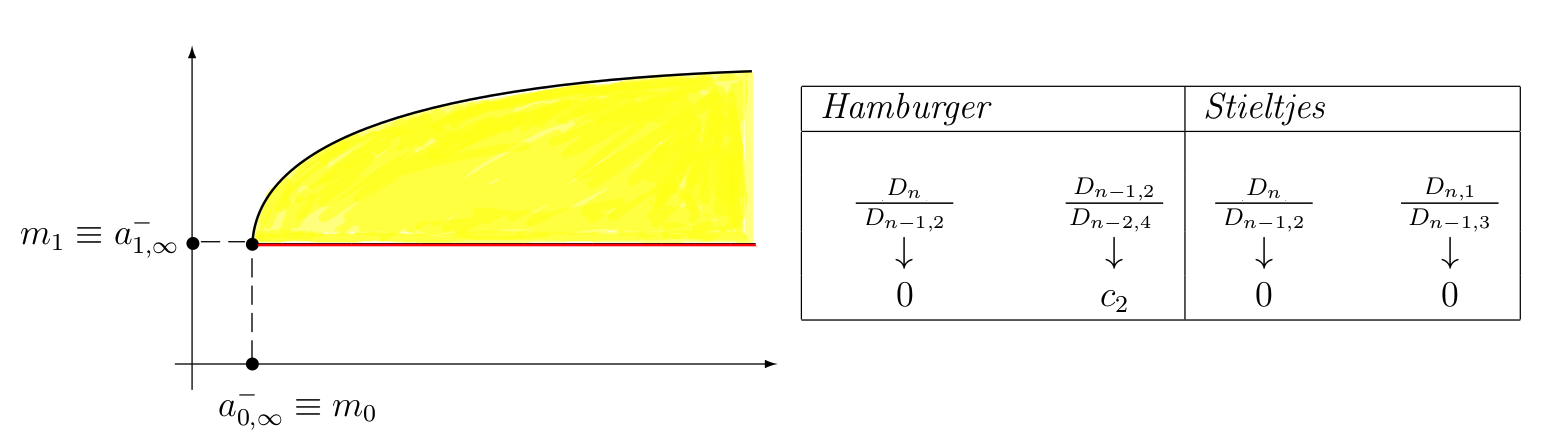}
		\caption{Parabolic region for S-det with H-det - (b)}
		\label{fig_6: H-Det-S-Det (b)_New (1)}
	\end{figure}

\noindent 
(c) and (d) \	 Here the limit parabolic regions are rays so that from the relation $a_{1,\infty}^{-}=m_1$ it 
follows that 	$\lim_{n \to \infty} \lambda_1(H_{n,1})=0$. Hence in both, (c) and (d), we have  that 
$a_{0,\infty}^{-}<m_0$. 
Graphically, see the red lines in Figure 3 
and in Figure 5, 
respectively.
It is interesting to mention that, if having S-determinate and H-determinate, each of the relations 
$\lim_{n\to\infty}\lambda_1(H_{n,1})=0$ and $\lim_{n\to\infty}\lambda_1(H_{n,1})=c_1$ may occur.

\smallskip
The findings in Cases 1 - 3 above can be summarized as follows:  

\smallskip
\noindent
$\bullet$ {\it If $\mu$ is H-indeterminate, equivalently, if \ $\lim_{n \to \infty}\lambda_1(H_{n})=c_0 > 0$, then the condition 
 $\lim_{n \to \infty}\lambda_1(H_{n,1})=0$ is  {\bf \emph{necessary and sufficient}} for $\mu$ to be S-determinate.}

\noindent
$\bullet$ {\it If $\mu$ is H-determinate, then $\lim_{n \to \infty}\lambda_1(H_{n})=0$ is a 
{\bf \emph{necessary and sufficient condition }} 
for $\mu$ to be S-determinate. Notice that each of the relations $\lim_{n\to\infty}\lambda_1(H_{n,1})=0$ may $\lim_{n\to\infty}\lambda_1(H_{n,1})=c_1>0$ may occur.}

\smallskip
The arguments used in Cases 1, 2, 3 above can be alternatively expressed in terms of the smallest eigenvalues of Hankel 
matrices thus arriving at a result which is equivalent to the known result in Berg-Thill (1991), Proposition 2.3.  

\smallskip\noindent
{\bf Theorem 6.} {\it	A Stieltjes moment sequence $\{m_k\}_{k=0}^{\infty}$ 
	corresponds to exactly one measure on the positive real axis {\bf \emph{if and only
	if}} the smallest eigenvalues of either $H_{n}$ or of $H_{n,1}$ tend to 0, as $n\to\infty$, i.e.
	\begin{equation}\label{SThm}
	\lim_{n\to\infty} \lambda_{1}(H_{n}) = 0 \quad or \quad      \lim_{n\to\infty} \lambda_{1}(H_{n,1}) = 0.
	\end{equation}
}

\smallskip
\noindent
{\bf Remark 2.} 
If we do not involve the lower bounds for the smallest eigenvalues of Hankel matrices,   
 Theorem 6 
can be easily proved  by combining the Krein-Nudelman's statement
used above, with the main result in Berg-Chen-Ismail (2002). 

Indeed, from Krein-Nudelman's result, we have that the measure $\mu$ with moments $\{m_{k}\}_{k=0}^{\infty}$
is S-indeterminate    {\it if and only if } \ H-indeterminate are both the measure $\mu$ with moments $\{m_{k}\}_{k=0}^{\infty}$ and 
the measure $\mu_1$ with the shifted moments $\{m_{k+1}\}_{k=0}^{\infty}$. 
The opposite statement is:  $\mu$ with  $\{m_{k}\}_{k=0}^{\infty}$ is S-determinate {\it if and only if } 
either $\mu$ with  $\{m_{k}\}_{k=0}^{\infty}$  is H-determinate, or $\mu_1$ with 
 $\{m_{j+1}\}_0^{\infty}$ is H-determinate. 
In terms of Theorem 1.1 in Berg-Chen-Ismail (2002), the last statement sounds as follows: 
$\mu$ with $\{m_{k}\}_{k=0}^{\infty}$ is S-determinate {\it if and only if} \  
$	\lim_{n\to\infty} \lambda_{1}(H_{n}) = 0$  or  $ \lim_{n\to\infty} \lambda_{1}(H_{n,1}) = 0$. $\Box$

\smallskip
We turn now to an important result relating S-determinacy and H-determinacy. Such a result is proved by 
Schm\"udgen (2017), Corollary 8.9, p. 183, in the framework of the operator-theoretic approach and by Heyde (1963), Theorem 1, p. 91, 
using continued fractions. We give a different short proof involving limit parabolic regions. 

\smallskip\noindent
{\bf Theorem 7.} {\it Suppose $\{m_{k}\}_{k=0}^{\infty}$ is a Stieltjes moment sequence associated with the measure $\mu$. 
	If $\mu$ is S-determinate  with zero mass at zero, $\mu(\{0\})=0$, then $\mu$ 
	considered on $\R$, is also H-determinate. 
}

\noindent
{\bf Proof.} Since $\mu$ is an S-determinate measure on $[0,\infty)$, there are two options. One is that $\mu$ is 
H-indeterminate. Then, referring to  Remark 1 after Case 2 above, $\mu$ must have a mass at 0, which is not the case. Thus it remains 
the second option for $\mu$, namely, that $\mu$ is H-determinate. Indeed, if turning to  
Theorems 1, 2, 3 and 4, we see that the limit parabolic regions in items (a), (b) and (d), see Case 3 above,  
are compatible  with the statement of Theorem 7. Notice, however, the parabolic region (no figure) in item  (c) has to be excluded, 
because of the appearance of a mass $\mu(\{0\})>0$, which contradicts the assumption.  
$\Box$

\bigskip
All measures/distributions satisfying the conditions of Theorem 7 are related to shifted Hankel matrices whose smallest eigenvalues 
may have different limits, as $n \to \infty,$ e.g.,  
 $\lambda_1(H_{n,1}) \to c_1 > 0$ or $\lambda_1(H_{n,1}) \to 0$   Clearly, Theorem 7 can be 
formulated in other equivalent forms. 

\smallskip
It is useful to provide here a result of Heyde (1963), his `Theorem B', which we paraphrase as follows: 

\smallskip\noindent
{\bf Theorem 8.} {\it 
  We are given a Stieltjes moment sequence  $\{m_0=1, m_{k}\}_{k=1}^{\infty} $ and let the associated 
 	probability measure $\mu$ be S-determinate with no mass at zero: $\mu(\{0\})= 0$. 
	
	Suppose that for fixed $\delta \in (0,1),$ a mass $\mu (\{0\})=\delta$ has been `added' at the origin $0$ and the distribution $\mu$ has been renormalized. 
	
	Then it is possible that the new moment sequence  
\[
	\{m_0^*=1, \ m_k^* := m_k/(1+\delta)\}_{k=1}^{\infty}
\]
generates a new distribution, say $\mu^*$, which is S-indeterminate. 
}

\noindent
{\bf Proof.} Indeed, from Theorem 7, the measure $\mu$ with moments $\{m_0=1, m_{k}\}_{k=1}^{\infty} $ 
is S-determinate and H-determinate and these properties are in agreement with the conclusions from the limit parabolic regions  in items 
(a), (b) and (d). The assumption $\mu (\{0\})>0$ changes the picture. The normalized measure $\mu^*$ for the new sequence	 
$\left\{m_0^{*}=1, m_k^*=\frac{m_{k}}{1 + \delta} \right\}_{k=1}^{\infty} $ is compatible with the following subcases: 
	(i) $\mu^*$ is H-indeterminate and S-indeterminate; (ii) $\mu^*$ is H-indeterminate and S-determinate;
	 (iii) $\mu^*$ is H-determinate and S-determinate. 
	 In fact, subcase  (i) proves Theorem B in Heyde (1963).
$\Box$

It is useful to add a few words. The statement in Heyde's Theorem B means that the shifted Hankel matrix $H^*_{n,1}$ based on the new moment 
sequence $\{m_k^*\}_{k=0}^{\infty}$ has a smallest eigenvalue $\lambda_1(H^*_{n,1})$ such that  $\lambda_1(H_{n,1}^*) \to c_1^* >0$. 
If we look at all subcases (i), (ii) and (iii), 
we see that, as $n \to \infty$,  either $\lambda_1(H^*_{n,1})\to 0$ or $\lambda_1(H^*_{n,1})\to c_1^*.$

\smallskip

\section{More on the lower bounds} 

Let us start with the Hamburger case. It is well-known and we have seen in the previous sections that there are some 
quantities which are important for  
deciding whether or not a measure is determinate or indeterminate. The modern technology allows `easily' to perform computations 
with good accuracy. Thus, in principle, having computed some quantities, would allow to make either preliminary or definite conclusions.  

Recall, writing below $H_n(m_0), n=1, 2, \ldots,$  means that $H_n$ are the basic Hankel matrices defined for the moment sequence $\{m_k\}_{k=0}^{\infty}.$ 
On several occasions we have introduced and used the numbers $b_k(m_0;n)$ for $k=1, \ldots, n$, and $b_{n+1}(m_0; n)$, see Section 3. 
It turns out, as $m_0 \to \infty$, their asymptotic behavior is different, 
namely: we have $b_k(m_0;n)\to 0$ for $k=1, \dots, n$, while $b_{n+1}(m_0; n) \to 1$. 
 As a consequence, if $m_0$ is finite, then for each $n$, it holds that $b_1(m_0; n)>0$ strictly. 

\smallskip\noindent
{\bf Proposition 1.} {\it 
Suppose that the following conditions are satisfied: 
\begin{equation}\label{boundm1}
  b_k(m_0; n) \leq  b_1(m_0; n) \quad \mbox{ for } \quad k=2, \ldots, n+1.
\end{equation}
Then for each fixed $n$, the lower bound \eqref{lamda1A} 
for the smallest eigenvalue $\lambda_1 (H_n(m_0))$ 
can be compared with the lower bound, {\rm BCI}, derived in Berg-Chen-Ismail (2002), eq. (1.14)-(1.15)). The relationship is as follows: 
\begin{equation} \label{BCI-lambda}
  {\rm BCI} := \Big(\frac{1}{2\pi } \int_0^{2\pi} \sum_{k=0}^{n} |p_k(e^{i\theta})|^2 d\theta\Big)^{-1} \leq 
	b_1(m_0; n)(m_0 - a_{0,n}^{-}) \leq  \lambda_1(H_{n}(m_0)).
\end{equation}
}

\noindent
{\bf Proof}. Indeed, from Berg-Chen-Ismail (2002), eqns. (1.12)-(1.13)), it follows  
	\begin{equation}
	\frac{1}{ \lambda_1(H_{n}(m_0)) }\leq \sum_{k=1}^{n+1}\frac{1}{ \lambda_k(H_{n}(m_0)) }= {\rm Tr}(H_{n}^{-1}(m_0) ) = \frac{1}{2\pi } \int_0^{2\pi} \sum_{k=0}^{n+1} |p_k(e^{i\theta})|^2 d\theta.
	\notag
	\end{equation}
To get the most left standing term we have neglected the quantity $\sum_{k=2}^{n+1} \frac{1}{\lambda_k(H_{n}(m_0))}$.
  From relation \eqref{Weyl2} we find that 
	\begin{equation}
	\frac{1}{ \lambda_1(H_{n}(m_0)) } \leq \frac{1}{b_1(m_0;n)) }\,\sum_{k=1}^{n+1}
	\frac{b_k(m_0;n))}{\lambda_k(H_{n}(m_0))} = \frac{1}{b_1(m_0;n)}\cdot \frac{1}{m_0- a_{0;n}^{-} }.  
	\notag
	\end{equation}
In this inequality neglected/omitted is the quantity $\frac{1}{b_1(m_0;n)) }\sum_{k=2}^{n+1} \frac{b_k(m_0;n)}{\lambda_k(H_{n}(m_0))}$. 
Conditions \eqref{boundm1} show that for the two neglected quantities we have  
\begin{equation}
	\frac{1}{b_1(m_0;n)}\,\sum_{k=2}^{n+1} \frac{b_k(m_0;n)}{\lambda_k(H_{n}(m_0))} \leq	\sum_{k=2}^{n+1} \frac{1}{\lambda_k(H_{n}(m_0))}.
	\notag
\end{equation}
Thus, the relationship (4.2) between the two lower bounds is established. $\Box$

Having  (4.2), we can take a limit, as $n \to \infty,$ by preserving all relations.  
We use the fact that the numbers $a_{0,n}^-$ have a limit, as $n \to \infty$, and we used for it the notation $a_{0,\infty}^-$.
The conclusion, after passing to the limit, is as follows: 
\begin{equation}\label{bounds}
		\Big(\frac{1}{2\pi } \int_0^{2\pi} \sum_{k=0}^{\infty} |p_k(e^{i\theta})|^2 d\theta\Big)^{-1} \leq b_1(m_0;\infty)(m_0 - a_{0,\infty}^{-}) \leq 
		 \lim_{n \to \infty} \lambda_1(H_{n}(m_0)). 
	\end{equation}

\smallskip
Relations \eqref{BCI-lambda} and \eqref{bounds} can eventually be validated through numerical examples. However first the assumptions \eqref{boundm1} have 
to be verified and guaranteed. For practical purposes, we can take $n$ large enough, and follow three steps:

\noindent
{\it Step 1.} The quantity $\rho_0$=:$\frac{1}{2\pi } \int_0^{2\pi} \sum_{k=0}^{\infty}|p_k(e^{i\theta})|^2 d\theta$ to be 
replaced by ${\rm Tr}\,(H_{n}^{-1} )$. 
 
\noindent
{\it Step 2.} The difference   $m_0 - a_{0,\infty}^{-}$ to be replaced by $1/H_{n}^{-1}(1,1).$ This is justified by the fact that 
    $m_0 - a_{0,n}^{-}=D_{n}/D_{n-1,2} = 1/H_{n}^{-1}(1,1)$.

\noindent
{\it Step 3.} The number $b_1(m_0;\infty)$ will also be replaced appropriately. The matrix $H_n$ is symmetric and diagonalizable and let us 
assume that $\lambda_1(H_{n})$ is a simple eigenvalue with its associated eigenvector $v=v({m_0;n})$. We use the notation $H_n(x)$ to indicate that the moment 
		$m_0 = x$ can `vary', so we find that $b_1(m_0;n)= 
		\frac{\rm d}{{\rm d}x}\lambda_1(H_{n}(x))|_{x=m_0}= \frac{v^T E v }{v^T v }$ 
   	 (see Golub-Van Loan (1996), p. 323). Then  for large $n$, 
   	$ b_1(m_0;\infty)$ can be replaced by $\frac{ v^T E v }{v^T v }.$ 
   	The eigenvector $v$ can be calculated efficiently by using the Inverse Power method.

\smallskip
A similar approach works also in the Stieltjes case. It is based on the statement by Krein-Nudelman, quoted before. The lower bound  \eqref{bounds} involves the matrices  $\{H_{n}(m_0)\}_{n=1}^{\infty}$, the moment sequence $\{m_0, m_1, m_2, \ldots \}$ and the smallest 
eigenvalues $\lambda_1 (H_n(m_0))$. Now, an analogue to \eqref{bounds}  
can be established for the shifted items $\{H_{n,1}(m_1)\}$, $\{m_1, m_2, \ldots \}$ and $\lambda_1 (H_{n,1}(m_1)).$ 
And,  we have to use $\{{\tilde p}_k(x)\}$, the sequence of orthonormal polynomials for the measure $\mu_1$, 
where ${\rm d}\mu_1 = x{\rm d}\mu$. 

\smallskip
\noindent
{\bf Remark 3.}
As we will see below, all numerically computed lower bounds in both Hamburger and Stieltjes cases look tight. 
In the H-indeterminate case, we have for $\lambda_1(H_n)$ a lower bound, a positive constant. 
This agrees with a result in Berg-Szwarc (2011), Theorem 4.4. These authors proved that, as $n\to\infty,$
the $k$th smallest eigenvalue of $H_n$, $\lambda_k(H_{n}(m_0))$, tends rapidly to infinity with $k$.  
This means that the two neglected terms $\sum_{k=2}^{n+1} \frac{1}{\lambda_k(H_{n}(m_0))}$
and $\sum_{k=2}^{n+1} \frac{ b_k(m_0;n)} {\lambda_k(H_{n}(m_0))}$, see the proof of Proposition 1, are indeed small.
In all cases, we have to be careful when using quantities such as 
${\rm Tr}(H_{\infty}^{-1})=\rho_0=:\frac{1}{2\pi } \int_0^{2\pi} \sum_{k=0}^{\infty} |p_k(e^{i\theta})|^2 {\rm d}\theta$ 
ensuring they are bounded.

\section{Numerical illustrations} 

While up to this place of our paper, we have used traditional terminology, notations and arguments from Analysis, now we turn to standard 
 probabilistic terminology and arguments. The main is the same, what is different is clear.

We have chosen two popular and frequently used probability distributions, namely, the Weibull distribution, which includes the exponential distribution, and 
the Lognormal distribution. Their (in)determinacy properties are well described and available in the literature. 
The reader can consult, e.g., 
Lin (2017) or Stoyanov-Lin-Kopanov (2020). Our goal now is to use the lower bounds for the 
smallest eigenvalues of Hankel matrices and, in a sense, confirm these (in)determinacy properties.    
 We give two examples. 
Example 1 is related to the content of the paper by Chen-Lawrence (1999) dealing with the weight function $w(x)=\exp (-x^{\beta}), \ x>0, \beta > 0.$ 
Notice, after normalizing $w$, it becomes the density function of the Weibull distribution (also called `generalized gamma distribution').   
Example 2 is similar to Example 3.1 in Berg-Chen-Ismail (2002), in which 
the authors start with the weight function $w(x)=x\,f(x), \ x>0$, where $f$ is the standard lognormal density. The treatment in these two 
papers is entirely analytic, no probabilistic notions involved.

\smallskip
\noindent
{\bf Example 1.}  (Weibull distribution). We say that a random variable $X$ has a {\it Weibull distribution} with parameter $\beta > 0,$  
$X \sim {\rm Wei}(\beta)$, if its probability density function is of the form 
\[
f(x)= c_{\beta}\,{\rm e}^{-x^{\beta}}, \ x > 0; \quad f(x)=0, \ x \leq 0.
\] 
Here $c_{\beta}$ is the normalizing constant. We easily see that all moments  
  $m_k := 
	\mathsf{E}[X^k] = \int_0^{\infty} x^k f(x)\,{\rm d}x, \ k=0, 1, 2, \ldots,$ are finite; $c_{\beta}$ and $m_k$ can be expressed 
	via the Euler gamma function. 	
	The (in)determinacy property of $X$ depends 
on the value of $\beta$. It turns out, $\beta = \frac12$ is the boundary point: if $\beta \geq \frac12$, $X$ and its distribution 
${\rm Wei}(\beta)$ are determinate, while they are indeterminate for any $\beta \in (0, \frac12).$ 

Consider now $\{m_k\}_{k=0}^{\infty}$ as a Stieltjes moment sequence and as a Hamburger moment sequence. 
We want to make the above conclusions by computing the lower bounds of the smallest eigenvalues of the Hankel matrices  $H_n$   and $H_{n,1}.$ 

As an illustration, assume that $\beta < \frac12$, expecting to obtain S-indeterminacy and also H-indeterminacy. These conclusions 
are correct if  based, e.g.,  on specific computations performed for $\beta = 0.45.$ Here are our conclusions: 

\smallskip
\noindent
$\bullet$ \ ${\rm Wei}(0.45)$ \ is \ H-{\it indeterminate}, which follows from the relations: 

$0.3323= \frac{1}{\rho_0}< 0.3402= b_1(m_0;\infty)(m_0-a_{0,\infty}^{-})< 0.3404=\lim_{n \to \infty}\lambda_1(H_{n})$. 

\noindent
$\bullet$ \ ${\rm Wei}(0.45)$ \ is \ S-{\it indeterminate}, since: 

$1.2588= \frac{1}{\rho_1}<  1.26175= b_1(m_1;\infty)(m_1- a_{1,\infty}^{-})<   1.26177=
\lim_{n \to \infty}\lambda_1(H_{n,1})$.
 
\vspace{0.2cm}
As a continuation, take $\beta =1,$ so we deal with a random variable $Y \sim {\rm Exp}(1),$ the exponential distribution 
with parameter 1, its density function is ${\rm e}^{-x}, \ x > 0.$ All moments of $Y$ are finite, $m_k={\mathsf E}[Y^k]=k!, \ k=1,2, \ldots.$  
In a few different ways we can show that $Y$, hence also ${\rm Exp}(1)$, is determinate. 

Moreover, for any power $Y^r, \ r >0,$ 
we easily find the density function, hence the distribution function (the measure), and see that all moments $m_k(Y^r), \ k=1, 2, \ldots,$ are finite; 
they are expressed via the Euler gamma function. The interesting property is that 
$Y^r$  is determinate for $r \in [0,2]$, and indeterminate for $r > 2.$ 

These conclusions can be derived from computed lower bounds of the smallest eigenvalues of Hankel matrices.     
We can write the matrices $H_n, \ H_{n,1}$ and compute that if $r \in (0,2]$, then  $\lambda_1(H_n) \to 0$ and $\lambda_1(H_{n,1}) \to 0$  
for large $n$.  This confirms that indeed $Y^r$ is S-determinate and also H-determinate. It is not surprising to observe that if $r$ is `close' 
to the boundary $r=2$, the convergence to zero is quite slow. 

It is instructive to make one step more by considering the random variable $Z$, where 
\[
Z = Y^3, \ \mbox{ for } \ Y \sim {\rm Exp}(1); \quad \mathcal{L}(Z) = {\rm Law}(Z). 
\] 
Notice, the number 3 is the smallest positive integer power such that 
$Y^3$ is indeterminate. Its moments are $m_k(Z) = \mathsf{E}[Z^k] = \mathsf{E}[Y^{3k}] = (3k)!, \ k=1, 2, \ldots$ 
The sequence $\{(3k)!\}_{k=0}^{\infty}$, being a Stieltjes moment sequence, can be considered also as a Hamburger moment sequence. We want to draw a 
conclusion for $\mathcal{L}(Z)$ based on calculated lower bounds of the smallest eigenvalues of the corresponding Hankel matrices. 
  With a reasonable accuracy of the computations, we arrive at the following conclusions:

 \smallskip\noindent
 $\bullet$ \ $\mathcal{L}(Z)$ \ is \ H-{\it indeterminate}, because  
 
 $  0.886774 = \frac{1}{\rho_0}< 0.8911283  = b_1(m_0;\infty)(m_0 - a_{0,\infty}^{-})<   0.8911307=\lim_{n \to \infty}\lambda_1(H_{n}).$
 
 \noindent
 $\bullet$ \ $\mathcal{L}(Z)$ \ is \ S-{\it indeterminate}, since    

 $ 3.00372= \frac{1}{\rho_1}<  3.003919   = b_1(m_1;\infty)(m_1- a_{1,\infty}^{-}) \leq  3.003919  = \lim_{n \to \infty}\lambda_1(H_{n,1}).$

\vspace{0.2cm}\noindent
 {\bf Example 2}.  (Lognormal distribution). We say that the random variable $\xi$ follows a {\it lognormal distribution}, \ $\xi \sim Log\,\mathcal{N}$, if 
its density function is 
\[
f(x) =\frac{1}{\sqrt{2\pi}}\cdot \frac{1}{x}\cdot \exp\left(- \frac12\,(\ln x)^2\right), \  x>0; \quad f(x) = 0, \ x \leq 0.
\]
All moments are finite, and  $m_k = \mathsf{E}[\xi^k] = {\rm e}^{k^2/2}, \ k=1, 2, \ldots$ Notice that   
$Log\,\mathcal{N}$ is the best known moment indeterminate absolutely continuous probability distribution. 

Let us draw the indeterminacy property from computed lower bounds for the smallest eigenvalues of the corresponding Hankel matrices. 
  So,  $\{{\rm e}^{k^2/2}\}_{k=0}^{\infty}$ being a Stieltjes moment sequence can also be considered as a Hamburger moment sequence. 
  With a reasonable computational accuracy, our results and conclusions are as follows:
 
\smallskip
 \noindent
 $\bullet$ \ $Log\,\mathcal{N}$ \ is \ H-{\it indeterminate}, which follows from the relations  
 
 $ 0.400108= \frac{1}{\rho_0}< 0.434605  = b_1(m_0;\infty)(m_0 - a_{0,\infty}^{-})<  0.441872 =\lim_{n \to \infty}\lambda_1(H_{n}).$
 
 \noindent
 $\bullet$ \ $Log\,\mathcal{N}$ \ is \ S-{\it indeterminate}, because

 $0.80338 = \frac{1}{\rho_1}< 0.817366  = b_1(m_1;\infty)(m_1- a_{1,\infty}^{-}) \leq 0.8176197  = \lim_{n \to \infty}\lambda_1(H_{n,1}).$

 \section{Brief concluding comments}

We exploit the geometric interpretation of indeterminacy conditions and fine properties of the eigenvalues of perturbed symmetric matrices 
and give a unified exposition of a bunch, or a series of classical results in the problem of moments. 
Our paper throws an additional light to the phenomena `uniqueness' and `nonuniqueness' (determinacy and indeterminacy) 
of  measures in terms of their moments. In general, it is always useful to have in our disposal  
different approaches, ideas and techniques which lead either to the same final conclusions or allow to establish new results.  
This  enhances the theory and makes it more applicable  
by providing the freedom to choose and use properly the most appropriate tools.

Despite the practical difficulty, in fact the impossibility, to check  the determinacy or indeterminacy conditions involving the smallest 
eigenvalues of Hankel matrices, 
the results discussed in this paper have been, are and will remain fundamental in mathematics. The numerical illustrations given above 
indicate that at least some of these 
results can be adopted and used in applications when dealing with specific measures/probability distributions.

\newpage\noindent
{\Large {\bf Acknowledgments}}

\vspace{0.3cm}
We are grateful to the editors and the technical staff for their understanding and support during the preparation of this paper. Our thanks  
are addressed to the two anonymous referees for their useful relevant comments.

\vspace{0.6cm}\noindent
{\bf \Large{References}}  

{\small
\vspace{0.1cm}
\begin{description}

\item
Akhiezer, N.I.  \textit{ The Classical Moment Problem and Some Related Questions in Analysis}. Oliver and Boyd: Edinburgh, 1965. (Original edition, Nauka, Moscow, 1961.)

\item
Berg, C.; Christensen, J.P.R.; Ressel, P. {\em Harmonic Analysis on Semigroups. Theory of Positive and Related Functions.} Springer: New York, 1984.  

\item 
Berg, C.; Thill, M. A density index for the Stieltjes moment problem. (Symposium on Orthogonal Polynomials, Erice, Sicily). In: Polynomials and Their 
Applications. Eds C. Brezinski, L. Gori, A. Ronvaux. Baltzer AG Sci. Publ. Co., IMACS, {\bf 1991}, pp. 185--188.   

\item
Berg, C.; Chen, Y.; Ismail, M.E.H.   
Small eigenvalues of large Hankel matrices: the indeterminate case. {\em  Math. Scand.} {\bf 2002}, \ {\it 91}, 67--81.

\item
Berg, C.; Christensen, J.P.R.  Density questions in the classical theory of moment.
{\em  Ann. Inst. Fourier} \ {\bf 1981}, {\it 31},  99--114.

\item
Berg, C.;  Szwarc, R.  The smallest eigenvalue of Hankel matrices.  
{\em Construct. Approx.} \ {\bf 2011}, {\it 34}, 107--133. 

\item
Berg, C.; Valent, G.  The Nevanlinna parametrization for some indeterminate Stieltjes moment problems associated 
with birth and death processes. {\em   Methods Appl. Anal.} \ {\bf 1994}, {\it 1},   169--209.

\item 
Chen, Y.; Lawrence, N.D.  Small eigenvalues of large Hankel matrices. {\em J. Phys. A} \ {\bf 1999},  {\it 32}, 7305--7315. 

\item
Chen, Y.; Sikorowski, J.; Zhu, M.  Small eigenvalues of large Hankel matrices at critical point: Comparing  conjecture 
of parallelised computation. {\em Appl. Math. Comput.} {\bf 2019}, {\it 363}, no. 124628. 

\item 
Chihara, T.S.  On indeterminate Hamburger moment problem. {\em Pacific J. Math.} \ {\bf 1968}, {\it 27}, 475--484.

\item
Golub, G.H.; Van Loan, C.F.   {\it Matrix Computations}. 
Johns Hopkins Univ. Press: Baltimore (MD), 1996.

\item 
Hamburger, H.  {\"U}ber  eine Ertweiterung des Stieltjesshen Momentenprobemes. {\it Math. Annalen} \ {\bf 1920}, 
{\it 81}, 235--317;   {\it 82}, 120--164, 168--187. 

\item
Heyde, C.C.   Some remarks on the moment problem. 
{\em  Quart. J. Math.}   (2) \ {\bf 1963}, {\it 14}, 91-–96.

\item 
Janssen, O.; Mirbabayi, M.; Zograf, P.  Gravity as an ensemble and the moment problem. {\it JHEP (SISSA)} (Springer) {\bf 2021}, 
{\it 06}:184, 19 pp. 

\item
Krein, M.G.; Nudelman, A.A.  {\it The Markov Moment Problem and Extremal Problems. (Ideas and Problems of 
P.L. Chebyshev and A.A. Markov and Their Further Development)}. Amer. Math. Soc.: Providence (RI), 1977.   
(Original edition, Nauka, Moscow, 1973.)

\item
Lin, G.D.  Recent developments on the moment problem. {\em J. Statist. Distrib. Appl.} \ {\bf 2017}, {\it 4}:1, 1--17. 

\item 
Lin, G.D.; Stoyanov, J.M. Moment Analysis of Probability Distributions: A Bunch of Old and New Checkable 
Determinacy Conditions. {\it Symmetry} {\bf 2023}, coming soon. 

\item 
Merkes E.P., Wetzel M.   A geometric characterization of indeterminate moment sequences. {\em Pacific J.  Math.} \ {\bf 1976}, 
{\it 65}, 409--419. 

\item 
Olteanu O. (2023). Symmetry and asymmetry in moment, functional equations and optimization problems. {\em Symmetry (MDPI)} {\bf 2023}, 
{\it 15}, 1471.  

\item
Schm\"udgen, K.  {\it The Moment Problem}. Graduate Texts in Mathematics {\bf 277}, Springer: Cham, 2017.

\item
Shohat, J.A.; Tamarkin, J.D.  {\it The Problem of Moments}. Math. Surveys No. {\bf 1}. Amer. Math. Soc.: Providence (RI), 1943. 

\item 
Sodin, S.  {\em Lecture Notes on the Moment Problem.} Department of Mathematics, Queen Mary University of London: London, 2019.\\
 Available at:  http://www.maths.qmul.ac.uk/$^{\sim}\hbox{s}_{-}$sodin/teaching/moment/clmp.pdf

\item 
Stieltjes, T.J.  Recherches sur les fractions continues. {\it Annales Fac. Sci. Univ. Toulouse} \ {\bf 1894}, {\it 8}, J1--J122; {\bf 1895}, {\it 9},  
A5--A47.   

\item
Stoyanov, J.; Lin, G.D.; Kopanov, P.  New checkable conditions for moment determinacy of probability distributions. 
{\em Theory Probab. Appl.} \ {\bf 2020}, {\it 65}, 497--509.

\item
Wilkinson, J.H.   {\it The Algebraic Eigenvalue Problem}. Clarendon Press: Oxford, 1985.

\item
Wulfsohn, A.  Measure convolution semigroups and noninfinitely divisible probability distributions. 
{\it J. Math. Sci. (NY)} \ {\bf 2005}, {\it 131}, 5682--5696. 

\end{description}

\end{document}